\documentclass[a4paper]{amsart}
\usepackage[all]{xy}

\usepackage{latexsym}
\usepackage{amsmath}
\usepackage{amssymb}

\newtheorem{theorem}{Theorem}
\newtheorem{lemma}{Lemma}
\newtheorem{proposition}{Proposition}
\newtheorem{corollary}{Corollary}

\theoremstyle{definition}

\newtheorem{remark}{Remark}

\newcommand{\p}{{\mathbb P}}
\newcommand{\g}{{\mathbb G}}
\newcommand{\df}{\operatorname{def}}

\begin{document}
\title{A remark on Zak's theorem on tangencies}
\author{Jos\'e Carlos Sierra}

\address {Instituto de Ciencias Matem\'aticas (ICMAT), Consejo Superior de
Investigaciones Cient\'{\i}ficas (CSIC), Campus de Cantoblanco,
28049 Madrid, Spain}

\email{jcsierra@icmat.es}

\thanks{Research supported by the ``Ram\'on y Cajal" contract RYC-2009-04999 and the project MTM2009-06964 of MICINN}

\begin{abstract}
We present a slightly different formulation of Zak's theorem on
tangencies as well as some applications. In particular, we obtain a
better bound on the dimension of the dual variety of a manifold and
we classify extremal and next-to-extremal cases when its secant
variety does not fill up the ambient projective space.
\end{abstract}
\maketitle

\section{Introduction}
In this note, we state Zak's theorem on tangencies (\cite[Theorem
0]{zak}, see also \cite[Ch.~I, Corollary 1.8]{zak2}) for
non-singular complex varieties in the following way:

\begin{theorem}[Reformulation of Zak's theorem on tangencies]\label{thm:main}
Let $X\subset\p^N$ be a non-degenerate manifold of dimension $n$,
and let $L\subset\p^N$ be a linear subspace of dimension $m$ which
is tangent to $X$ along a closed subvariety $Y\subset X$ of
dimension $r$. Then: $$r\leq\min\{m-n,\, \dim SX-1-n\}$$
\end{theorem}

The former statement appears to have some advantages. The new
inequality $r\leq\dim SX-1-n$ is vacuous if the secant variety of
$X$ (denoted by $SX$) fills up $\p^N$, but it is significant when
$SX\neq\p^N$. For instance, this is always the case if $m=N-1$. In
this setting, Zak's theorem on tangencies has several consequences
(see \cite{laz} for an account) that can be sharpened thanks to
Theorem \ref{thm:main}. Let $s:=\dim SX$ and $c:=N-s$.

\begin{corollary}\label{cor:main}
Let $X\subset\p^N$ be a non-degenerate manifold of dimension $n$,
and let $X^*\subset\p^{N*}$ denote its dual variety (of dimension
$n^*$). The following holds:
\begin{enumerate}
\item[(i)] The twisted normal bundle $N_{X/\p^N}(-1)$ is $k$-ample for $k\geq s-1-n$ (cf. \cite[Example
6.3.7]{laz}).
\item[(ii)] $n^*\geq n+c$ (cf. \cite[Corollary
3.4.20]{laz}). In particular, if $SX\neq\p^N$ then $X^*$ is a
singular variety.
\item[(iii)] If $s\leq 2n-1$ (resp. $2n-2$) then every hyperplane section of
$X$ is reduced (resp. normal) (cf. \cite[Corollary 3.4.19]{laz}).
\end{enumerate}
\end{corollary}

Some examples of manifolds satisfying the equality $n^*=n$ (and
hence $SX=\p^N$) are given by hypersurfaces in $\p^{n+1}$, Segre
embeddings $\p^1\times\p^{n-1}\subset\p^{2n-1}$, the Grassmannian
$\g(1,4)\subset\p^9$ and the $10$-dimensional spinor variety
$S_4\subset\p^{15}$. Moreover, these are the only examples under the
additional assumption $3n\leq 2N$ by \cite[Theorem 4.5]{ein} (cf.
Remark \ref{rem:ein}). On the other hand, if $SX\neq\p^N$ we show
that furthermore $n^*\geq n+c+1$ and manifolds satisfying the
equality are classified, giving a new characterization of the
Veronese surface:

\begin{theorem}\label{thm:veronese}
Let $X\subset\p^N$ be a non-degenerate manifold of dimension $n$. If
$SX\neq\p^N$ then $n^*\geq n+c+1$, with equality if and only if $X$
is either a curve or the Veronese surface in $\p^5$.
\end{theorem}

Going one step further, let us consider the next-to-extremal case
when $SX\neq\p^N$. If $n\leq 3$ it is easy to see (cf. Remark
\ref{rem:n^*}) that $n^*=n+c+2$ if and only if $X$ is either a
surface, or a dual defective threefold (i.e. a scroll over a curve),
or a secant defective threefold (see \cite{fuj2} for the
classification). On the other hand, for $n\geq 4$ we get the
following:

\begin{theorem}\label{thm:scroll}
Let $X\subset\p^N$ be a non-degenerate manifold of dimension $n\geq
4$. If $SX\neq\p^N$ and $n^*=n+c+2$, then $X$ is a scroll over a
manifold $W$ and $\dim W\leq 2$.
\end{theorem}

If $\dim W=2$ we will prove in Theorem \ref{thm:segre} that the
Segre embedding of $\p^2\times\p^{n-2}$ is the only scroll as in
Theorem \ref{thm:scroll}, so we actually get the following
refinement:

\begin{theorem}\label{thm:refinement}
Let $X\subset\p^N$ be a non-degenerate manifold of dimension $n\geq
4$ such that $SX\neq\p^N$. Then $n^*=n+c+2$ if and only if $X$ is
either a scroll over a curve, or (an isomorphic projection of) the
Segre embedding $\p^2\times\p^{n-2}\subset\p^{3n-4}$.
\end{theorem}

\section{Proofs}\label{sec:proof}
Theorem \ref{thm:main} is a consequence of the following application
of the Fulton-Hansen connectedness theorem \cite{f-h}. First, we
recall the definition of the relative tangent (resp. secant)
variety. Given a subvariety $Y\subset X$, we define
$T(Y,X):=\cup_{y\in Y}T_yX$, where $T_yX\subset\p^N$ denotes the
embedded tangent space to $X$ at $y\in Y$, and $S(Y,X)\subset\p^N$
as the closure of $\{z\in\p^N\mid \exists y\in Y, \exists x\in X
\,\text{with}\,\, z\in\langle y,x\rangle\}$.

\begin{theorem}\label{thm:f-h}
Let $X\subset\p^N$ be a non-degenerate manifold of dimension $n$,
and let $Y\subset X$ be a closed subvariety of dimension $r$. Then
either $\dim T(Y,X)=r+n$ and $\dim S(Y,X)=r+n+1$, or else
$T(Y,X)=S(Y,X)$.
\end{theorem}

\begin{proof}
See \cite[Ch.~I, Theorem 1.4]{zak2}.
\end{proof}

We can now prove our results:

\begin{proof}[Proof of Theorem \ref{thm:main}]
Let $L\subset\p^N$ be a linear subspace of dimension $m$ which is
tangent to $X$ along $Y$. Then $T(Y,X)\subset L$, but
$S(Y,X)\not\subset L$ as $X\subset\p^N$ is non-degenerate. Therefore
$T(Y,X)\neq S(Y,X)$, and hence $r+n=\dim T(Y,X)\leq \dim L=m$ by
Theorem \ref{thm:f-h}. But Theorem \ref{thm:f-h} also yields
$r+n+1=\dim S(Y,X)\leq\dim SX$, so $r\leq\min\{m-n,\, \dim
SX-1-n\}$.
\end{proof}

\begin{remark}
If $SX\neq\p^N$ and $m=N-1$, the case we are more interested in, the
new bound $r\leq s-1-n$ is sharp. For example, equality holds for
Severi varieties \cite[Ch.~IV]{zak2} when $Y\subset X$ is an
$n/2$-dimensional quadric.
\end{remark}

\begin{proof}[Proof of Corollary \ref{cor:main}]
According to Theorem \ref{thm:main}, the dimension of the fibres of
the second projection of the conormal variety
$\mathcal{P}_X:=\{(x,H)\mid T_xX\subset H\}\subset X\times\p^{N^*}$
is bounded by $s-1-n$. So $N_{X/\p^N}(-1)$ is $k$-ample for $k\geq
s-1-n$. This proves (i). Since $\dim\mathcal{P}_X=N-1$, we get
$n^*\geq n+c$. Assume now that $X^*$ is smooth. Then
$\dim(X^*)^*\geq n^*\geq n+c$, and Segre's reflexivity theorem
$(X^*)^*=X$ (\cite{seg}, see also \cite{kle} for a detailed account)
yields $c=0$, whence $SX=\p^N$ proving (ii). Part (iii) is an
immediate consequence of Theorem \ref{thm:main}.
\end{proof}

The main ingredients of the proof of Theorem \ref{thm:veronese} are
Zak's classification of Severi varieties and Ein's bound on the
defect of subcanonical manifolds \cite[Theorem 4.4]{ein}:

\begin{proof}[Proof of Theorem \ref{thm:veronese}]
Assume $n\geq 2$ and $n^*\leq n+c+1$. Let $\df(X)$ and $\delta(X)$
denote the dual and secant defect of $X\subset\p^N$, respectively.
As $\df(X):=N-1-n^*$ and $\delta(X):=2n+1-s$, the inequality
$n^*\leq n+c+1$ is equivalent to the inequality
$\df(X)+\delta(X)\geq n-1$. Since we assume $SX\neq\p^N$, we get
$\delta(X)\leq n/2$ by Zak's theorem on linear normality
\cite[Ch.~II, Corollary 2.11]{zak2} and equality holds if and only
if $X\subset\p^N$ is a Severi variety. Assume first that the Picard
group of $X$ is cyclic. Then $\df(X)\leq (n-2)/2$ by \cite[Theorem
4.4]{ein}, and hence
$$n-1\leq\df(X)+\delta(X)\leq\frac{n-2}{2}+\frac{n}{2}=n-1$$
implies that $X\subset\p^N$ is a Severi variety with
$\df(X)=(n-2)/2$. So $X$ is the Veronese surface in $\p^5$. Assume
now that the Picard group of $X$ is not cyclic. Then $\delta(X)\leq
2$ by the Barth-Larsen theorem, as otherwise $X\subset\p^N$ could be
isomorphically projected into $\p^{2n-2}$ (see for instance
\cite[Corollary 3.2.3]{laz}). Therefore $\df(X)\geq n-3$, whence
$\df(X)=n-2$ by Landman's parity theorem (unpublished, see
\cite[Theorem 2.4]{ein}) and $X\subset\p^N$ is a scroll over a curve
by \cite[Theorem 3.2]{ein}. This yields $1\leq\delta(X)\leq 2$,
contradicting Lemma \ref{lem:segre}.
\end{proof}

\begin{remark}\label{rem:n^*}
(i) The bound $n^*\geq n+c$ given in Theorem \ref{thm:main} is
equivalent to the bound $\df(X)+\delta(X)\leq n$. Furthermore, if
$SX\neq\p^N$ then $\df(X)+\delta(X)\leq n-1$ by Theorem
\ref{thm:veronese}. To the best of the author's knowledge, these
relations involving both dual and secant defects appear to be new.

(ii) For $n\geq 3$ the bound obtained in Theorem \ref{thm:veronese}
is equivalent to the bound $\df(X)+\delta(X)\leq n-2$. This can be
seen as a refinement of the inequality $\df(X)\leq n-2$ of Landman
and Zak when $SX\neq\p^N$ (cf. \cite[Ch.~I, Remark 2.7]{zak2}).

(iii) Besides curves and the Veronese surface, the bound $n^*\geq
n+c+2$ (or equivalently $\df(X)+\delta(X)\leq n-2$) is sharp.
Equality holds for surfaces with $\delta(X)=0$, threefolds with
$\delta(X)=1$, scrolls over curves with $\delta(X)=0$ and the Segre
embeddings $\p^2\times\p^{n-2}\subset\p^{3n-4}$ with $n\geq 4$. We
will prove in the sequel that these are actually the only ones.
\end{remark}

The key of the proof of Theorem \ref{thm:scroll} is a recent
characterization of scrolls among dual defective manifolds obtained
by Ionescu and Russo in \cite{i-r}:

\begin{proof}[Proof of Theorem \ref{thm:scroll}]
Let $n^*=n+c+2$, that is, $\df(X)+\delta(X)=n-2$. If $X$ is a scroll
over a manifold $W$ then $\delta(X)\leq 2$ by the Barth-Larsen
theorem. Thus $n-4\leq\df(X)=n-2\dim W$, so $\dim W\leq 2$. If $X$
is not a scroll, we can assume $\df(X)\leq(n+1)/3$ by
\cite[Corollary 3.7]{i-r} and $\delta(X)\leq(n-1)/2$ by the
classification of Severi varieties. Therefore,
$$n-2=\df(X)+\delta(X)\leq\frac{n+1}{3}+\frac{n-1}{2}$$ yields $n\leq
11$ and, in view of Landman's parity theorem, we get
$(n,\df(X),\delta(X))\in\{(5,1,2),(6,2,2),(9,3,4)\}$. The first two
cases are excluded by \cite[Theorems 5.1 and 5.2]{ein2}. In the
third case, $X$ is a Fano manifold of dimension $9$ with cyclic
Picard group generated by the hyperplane section and index
$(n+\df(X)+2)/2=7$ (see \cite[Lemma 4.2]{ein}), so it is ruled out
by Mukai's classification of Fano manifolds of coindex 3 (see
\cite{muk}).
\end{proof}

\begin{remark}\label{rem:ein}
(i) In a similar way, we can prove that a non-degenerate $n$-fold
$X\subset\p^N$ with $n^*=n$ is either the Segre embedding
$\p^1\times\p^{n-1}\subset\p^{2n-1}$, or else $4n+5\geq 3N$ and
equality holds if and only if $X$ is the $10$-dimensional spinor
variety $S_4\subset\p^{15}$. Since $n^*=n$ we deduce $SX=\p^N$ by
Corollary \ref{cor:main}, and hence $\delta(X)=2n+1-N$. We point out
that $n^*=n$ is equivalent to $\df(X)+\delta(X)=n$. If $X$ is a
scroll then $\delta(X)\leq 2$, whence $\df(X)\geq n-2$. Therefore
$\df(X)=n-2$, $\delta(X)=2$ and $X$ is the Segre embedding
$\p^1\times\p^{n-1}\subset\p^{2n-1}$ by Proposition
\ref{prop:scroll}. On the other hand, if $X$ is not a scroll then
$\df(X)\leq(n+2)/3$ and equality holds if and only if $X$ is the
$10$-dimensional spinor variety by \cite[Corollary 3.7]{i-r}. Thus
$\delta(X)\geq (2n-2)/3$, and hence $4n+5\geq 3N$, with equality if
and only if $X$ is the $10$-dimensional spinor variety
$S_4\subset\p^{15}$.

(ii) If $n^*=n$ and one furthermore assumes that $3n\leq 2N$ (cf.
\cite[Theorem 4.5]{ein}) then one also gets hypersurfaces and the
Grassmannian $\g(1,4)\subset\p^{9}$, as in \cite[Corollary
3.9]{i-r}.
\end{remark}

\begin{proof}[Proof of Theorem \ref{thm:refinement}]
In Theorem \ref{thm:scroll}, if $\dim W=2$ then $\df(X)=n-4$ and
hence $\delta(X)=2$, since $\df(X)+\delta(X)=n-2$. So we conclude in
view of Theorem \ref{thm:segre}.
\end{proof}

\section{A result on secant defective scrolls over surfaces}
In this section we prove the results on scrolls quoted in Section
\ref{sec:proof}. We say that $X_W\subset\p^N$ (or simply $X$) is a
\emph{scroll} if there exists a vector bundle $\mathcal E$ over a
manifold $W$ such that $X_W\cong\p_W({\mathcal E})$ and the fibres
of the map $\pi:X_W\to W$, that we denote by $F_w$ for $w\in W$, are
linearly embedded in $\p^N$. An equivalent definition of scroll is
the following. Let $\g(k,N)$ denote the Grassmannian of $k$-planes
in $\p^N$. Consider the incidence correspondence
$\mathcal{U}:=\{(\p^k,p)\mid p\in\p^k\}$ with projection maps
$\pi_1:\mathcal{U}\to\g(k,N)$ and $\pi_2:\mathcal{U}\to\p^N$. For
every subvariety $W\subset\g(k,N)$, we denote
$\mathcal{U}_W:=\pi_1^{-1}(W)$ and $X_W:=\pi_2(\mathcal{U}_W)$. Then
$X_W\subset\p^N$ is a scroll if and only if $W$ is smooth and
$\pi_2:\mathcal{U}_W\to X_W$ is an isomorphism. The following
consequence of Terracini's lemma \cite{ter} will be useful. Let
$\Sigma_z\subset X$ denote the \emph{entry locus} of $z\in SX$, that
is, the closure of the set $\{x\in X\mid\exists x'\in
X\,\text{with}\,\, z\in\langle x,x'\rangle\}$. We recall that
$\dim(\Sigma_z)=\delta(X)$ for general $z\in SX$.

\begin{lemma}\label{lem:useful}
Let $X\subset\p^N$ be a non-degenerate scroll over $W$ and let $z\in
SX$ be a smooth point. If $\Sigma_z\cap F_w\neq\emptyset$ for every
$w\in W$ then $SX=\p^N$.
\end{lemma}

\begin{proof}
Let $T_zSX\subset\p^N$ be the embedded tangent space to $SX$ at $z$.
For every $w\in W$ there exists some $x\in\Sigma_z\cap F_w$, so
$F_w\subset T_xX$. Then it follows from Terracini's lemma that
$X=\cup_{w\in W}F_w\subset\cup_{x\in\Sigma_z}T_xX\subset T_zSX$.
Since $X\subset\p^N$ is non-degenerate we deduce $T_zSX=\p^N$, and
hence $SX=\p^N$.
\end{proof}

The following lemma is well known. We include a short proof based on
Lemma \ref{lem:useful}:

\begin{lemma}\label{lem:segre}
Let $X\subset\p^N$ be a non-degenerate scroll over a curve. If
$\delta(X)>0$ then $SX=\p^N$.
\end{lemma}

\begin{proof}
Since $\dim\Sigma_z=\delta(X)>0$ for general $z\in SX$ and $\dim
W=1$, we deduce that $\Sigma_z\cap F_w\neq\emptyset$ for every $w\in
W$. Therefore $SX=\p^N$ by Lemma \ref{lem:useful}.
\end{proof}

Let $X_W\subset\p^N$ be an $n$-dimensional scroll. It follows from
the Barth-Larsen theorem that $\delta(X)\leq 2$. From now on, we
will focus on the extremal case $\delta(X)=2$. On the one hand, if
$W$ is a curve then $X_W\subset\p^N$ is the Segre embedding
$\p^1\times\p^{n-1}\subset\p^{2n-1}$ (see \cite[pp. 307--308]{kl}).
We prove this result in a more geometric and elementary way. The
idea of the proof is essentially due to Fyodor Zak:

\begin{proposition}\label{prop:scroll}
The only $n$-dimensional scroll over a curve in $\p^{2n-1}$ is the
Segre embedding $\p^1\times\p^{n-1}\subset\p^{2n-1}$.
\end{proposition}

\begin{proof}
For every $w\in W$, let
$\sigma_w:=\{g\in\g(n-1,2n-1)\mid\p_g^{n-1}\cap F_w\neq\emptyset\}$.
Then $\sigma_w$ is a hyperplane section of $\g(n-1,2n-1)$, embedded
by Pl\"ucker, and $w\in\sigma_w$ is a singular point of multiplicity
$n$. The intersection of $W$ and $\sigma_w$ is supported at $w$
since $\pi_2:\mathcal{U}_W\to X_W$ is injective. Moreover, $W$ and
$\sigma_w$ meet transversally at $w$ since $d\pi_2$ is injective. So
we deduce that $\deg X_W=\deg W=m_w(\sigma_w)\cdot m_w(W)=n$, where
$m_w(\sigma_w)$ and $m_w(W)$ denote the multiplicity of $\sigma_w$
and $W$ at $w$, respectively (see \cite[Corollary 12.4]{ful}).
Therefore, for every $w\in W$ there exists a hyperplane section
$\sigma_w$ in the Pl\"ucker embedding of $W$ such that the
intersection product $\sigma_w\cdot W=nw$. This property
characterizes the rational normal curve of degree $n$, so
$X_W\subset\p^{2n-1}$ is a non-degenerate (otherwise $F_w\cap
F_{w'}\neq\emptyset$ for every $w'\in W$) rational normal scroll of
degree $n$. Consequently, $X_W$ is the Segre embedding
$\p^1\times\p^{n-1}\subset\p^{2n-1}$.
\end{proof}

\begin{remark}\label{rem:segre}
The hypothesis of Proposition \ref{prop:scroll} can be weakened.
Arguing with a general $w\in W$, the same proof works if
$W\subset\g(n-1,2n-1)$ is an integral curve and
$\pi_2:\mathcal{U}_W\to X_W$ is an isomorphism (or even if
$\pi_2:\mathcal{U}_W\to X_W$ has finitely many double points).
\end{remark}

On the other hand, if $W$ is a surface there exists a complete
classification of scrolls with $\delta(X)=2$ only for $n=3$ (see
\cite{ott} and \cite[Proposition 4]{i-t}). We now prove that (an
isomorphic projection of) the Segre embedding
$\p^2\times\p^{n-2}\subset\p^{3n-4}$ with $n\geq 4$ is the only
scroll over a surface whose secant variety does not fill up the
ambient space. The main idea of the proof is to show that
$X\subset\p^N$ is swept out by a $2$-dimensional family of Segre
embeddings $\p^1\times\p^{n-2}$. More precisely, we prove that any
two fibres of the scroll, $F_w$ and $F_{w'}$, determine a Segre
embedding $\p^1\times\p^{n-2}$ in the linear span $\langle
F_w,F_{w'}\rangle=:\p_{ww'}^{2n-3}$ that they define.

\begin{theorem}\label{thm:segre}
Let $X\subset\p^N$ be a non-degenerate scroll of dimension $n$ over
a surface. If $\delta(X)=2$ and $SX\neq\p^N$ then $X$ is (an
isomorphic projection of) the Segre embedding
$\p^2\times\p^{n-2}\subset\p^{3n-4}$.
\end{theorem}

\begin{proof}
We claim that $\dim S(F_w,X)=2n-2$ for general $w\in W$. Since
$\langle F_w,F_{w'}\rangle=S(F_w,F_{w'})\subsetneq S(F_w,X)$ for
every $w'\in W$, we get $\dim S(F_w,X)\geq 2n-2$. Fix a general
$z\in SX$. We deduce from Lemma \ref{lem:useful} that $\Sigma_z\cap
F_w=\emptyset$ for general $w\in W$, and hence $z\notin S(F_w,X)$.
Therefore $S(F_w,X)\subsetneq SX$ proving the claim, as $\dim
SX=2n-1$. For every $w\in W$, consider the subvariety
$\mathcal{G}_w:=\{\langle F_w,F_{w'}\rangle\mid w'\in
W\}\subset\g(2n-3,N)$. If $\dim\mathcal{G}_w=2$ for general $w\in W$
then $S(F_w,X)\subset\p^N$ is a $(2n-2)$-dimensional subvariety
swept out by a $2$-dimensional family of $(2n-3)$-dimensional linear
subspaces, so $S(F_w,X)\subsetneq\p^N$ itself is a linear subspace.
This contradicts the non-degeneracy of $X\subset\p^N$. Thus
$\dim\mathcal{G}_w=1$ for general (and hence every) $w\in W$. In
particular, for every $w,w'\in W$ there exists an integral curve
$T_{ww'}\subset W$ such that $\langle F_w,F_{w'}\rangle=\langle
F_w,F_{w''}\rangle$ for every $w''\in T_{ww'}$. So
$X_{T_{ww'}}=\p^1\times\p^{n-2}\subset\p_{ww'}^{2n-3}$ by Remark
\ref{rem:segre}. Consequently, for every $w\in W$ and every $x\in
F_w$ there exists a $1$-dimensional family of lines each of them
meeting $F_w$ at $x$ and giving a $2$-dimensional cone $C_x\subset
X$. Since $F_w$ and $C_x$ are contained in $T_xX=\p^n$ we deduce
$\deg(C_x)=C_x\cdot F_w$. We claim that $C_x\cdot F_{w'}=1$ for
every $w'\in W$, and hence $\deg(C_x)=C_x\cdot F_w=C_x\cdot
F_{w'}=1$. Let us prove the claim. If $C_x\cdot F_{w'}\geq 2$ then
$T_xX=\langle F_w, C_x\cap F_{w'}\rangle\subset\langle
F_w,F_{w'}\rangle$. Therefore $T(F_w,X)\subset\langle
F_w,F_{w'}\rangle$, contradicting Theorem \ref{thm:main}. Since
$\deg(C_x)=C_x\cdot F_w=1$, we deduce that $C_x=\p^2$ for every
$x\in F_w$ and that $C_x$ is a section of $\pi:X_W\to W$ (in
particular, $W\cong\p^2$). Thus $X$ is a scroll over $\p^2$ and
$\p^{2n-2}$, respectively. So $X\subset\p^N$ is an isomorphic
projection of the Segre embedding
$\p^2\times\p^{n-2}\subset\p^{3n-4}$.
\end{proof}

\section*{Acknowledgements}
The author is grateful to Fyodor Zak for helpful comments and
encouragement.

\bibliography{bibfile}
\bibliographystyle{amsplain}

\end{document}